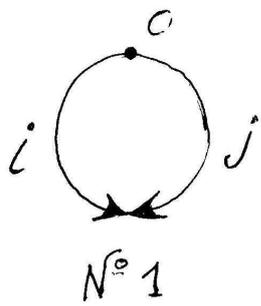

№ 1

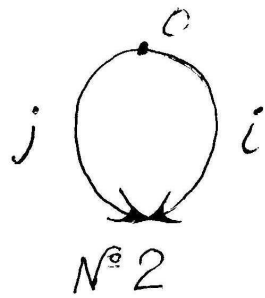

№ 2

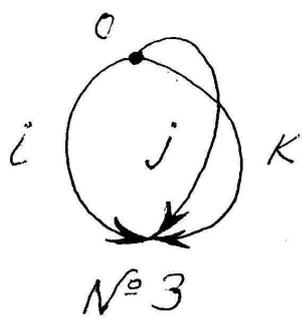

№ 3

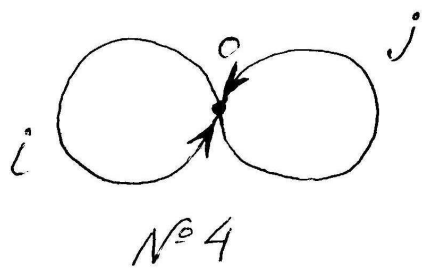

№ 4

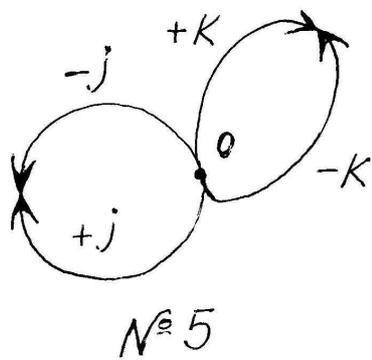

№ 5

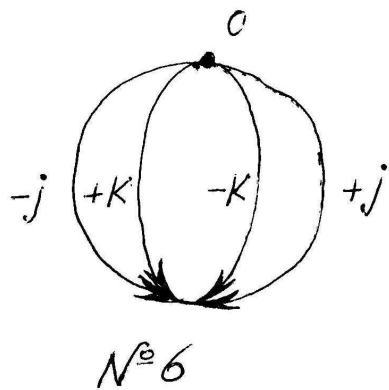

№ 6

И.В. Баяк

# О жесткой алгебраической структуре
# евклидовых пространств

***Резюме:*** Развит алгебраический формализм, в рамках которого показано, что группы автоморфизмов евклидовых пространств, сохраняющих ориентацию базиса, изоморфны группам топологических автоморфизмов соответственно факторизованных арифметических пространств.

### 0. Геометрия и алгебра компактных пространств

Поскольку представленный далее алгебраический формализм является базисным не только по отношению к геометрии евклидовых пространств, но и по отношению к так называемой геометрии компактных пространств, то укажем здесь первоисточник такого алгебраически геометрического соответствия, а именно, не ассоциативное по сложению «поле» $K = \langle \left]-2,2\right[; \oplus, \otimes, 0, 1 \rangle$, где $\varphi_1 \oplus \varphi_2 = |\varphi_1 + \varphi_2| \mod 2$, $\varphi_1 \otimes \varphi_2 = |\varphi_1 \cdot \varphi_2| \mod 2$, $|x| \mod 2 : R \to \pm S^1 \left(e^{i\pi\varphi} = e^{i\pi x}\right)$, которое реализуется на окружности и при этом имеет ограничения в области целых и натуральных чисел $K|_Z = \langle \{-1, 0, 1\}; \oplus, \cdot, 0, 1 \rangle$ и $K|_N = \langle \{0, 1\}; \oplus, \cdot, 0, 1 \rangle$, совпадающие с базовыми структурами алгебраического формализма, развиваемого в пункте 3, так что $B \approx K|_Z$ а $D \approx K|_N$ соответственно. Таким образом, декартово произведение $n$ экземпляров $K$, следует считать $n$-мерным компактифицированным арифметическим пространством и интерпретировать его как базовую структуру геометрии компактных пространств, воображаемую взглядом из 3-мернорго евклидова пространства как пучок из $n$ ориентированных окружностей, пересекающихся в одной точке.

### 1. Группа стрелочных подстановок и группы автоморфизмов некоторых групповых алгебр

Пусть $I = \{1, .., n\}$, $A = \{\pm 1, .., \pm n\}$, а $p: I \to A$ - такая инъекция, что $|p|: I \to I$ - биекция, причем $p'p: I \to A : p'p(i) = \text{sign} p(i) \cdot p'(|p(i)|)$, где $\text{sign} : A \to \{\pm 1\}$, $|\ |: A \to I$; тогда множество всех биективных по модулю инъекций $P = \{p\}$ относительно определенной ранее композиции, т.е. произведения $P \times P \to P : (p'p) \to p'p$, составляет группу стрелочных подстановок $n$-ой степени $P_n$, которая изоморфна группе переходных матриц порядка $n$, т.е. квадратных матриц с одним ненулевым элементом (+1 или -1) в каждой строке и каждом столбце. Добавим к этому, что стрелочная подстановка называется стрелочной транспозицией если это элементарная перестановка или элементарная инверсия, а четностью стрелочной подстановки $p$ называется четность композиции стрелочных транспозиций для перехода от тождественной подстановки к $p$. Таким образом, множество всех четных стрелочных подстановок $P_n^+$ составляет подгруппу, которая изоморфна подгруппе переходных матриц с единичным детерминантом. Вместе с тем, множество всех подстановок с четным числом инверсий $P_n^-$ составляет подгруппу, которая изоморфна подгруппе переходных матриц с четным числом отрицательных элементов. Пусть теперь семейство $I$ разбито на $m$ непересекающихся подсемейств, т.е. $I_J = \{I_j\}_J$, где $J = \{1,...,m\}, m < n$, $\bigcup_m I_j = I$, $\bigcap_m I_j = \varnothing$, $\text{card}(I_j) = n_j$, $\sum_m n_j = n$, тогда всякая простая и стрелочная подстановка может быть представлена как составная подстановка, состоящая из $m$ размещений по местам $I_j$, причем всякое стрелочное размещение по месту $I_j$ получено композицией элементарных размещений (элементарных перестановок и инверсий в подмножестве $I_j$ и элементарных замещений из подмножества $I \setminus I_j$), а четность размещения определяется четностью этой композиции. Аналогичное понятие четности простого размещения (отличающееся от стрелочного лишь отсутствием инверсий) позволяет определить детерминант подмножества строк квадратной матрицы как сумму произведений элементов, индексированных размещением столбцов по строкам и умноженных на знак четности этого размещения. В свою очередь, произведение знаков четностей всех $m$ размещений определяет четность составной подстановки, так что множество всех четно-составных стрелочных подстановок $P_n^\pm$ составляет подгруппу, которая изоморфна подгруппе переходных матриц с единичным произведением детерминантов $m$ прямоугольных матриц, соответствующих разбиению $n$ строк на $m$ подсемейств $I_j$.

Далее пусть в рамках процедуры порождения вышеуказанных групп и связанных с ними алгебраических структур принято, что $l = \begin{pmatrix} 0 & 1 \\ -1 & 0 \end{pmatrix}$, $t = \begin{pmatrix} 0 & 1 \\ 1 & 0 \end{pmatrix}$, тогда $P_2^+ = F(l)$, $P_2^- = F(\pm t)$, где $F(S)$ - свободная группа, образованная множеством $S$. Кроме того, генераторы $l$, $t$ формируют свободные алгебры над $R$, а именно, групповые алгебры $RP_2^+ = \begin{pmatrix} x & y \\ -y & x \end{pmatrix}$ и $RP_2^- = \begin{pmatrix} x & y \\ y & x \end{pmatrix}$, где $x, y \in R$, которым соответствуют группы их автоморфизмов, а именно, $L \equiv Aut\, RP_2^+ = GP_2^+ / R^+ = SO(2)$ и $T \equiv Aut\, RP_2^- = GP_2^- / R^+ = SO(1,1)$, где $GP_2^+$, $GP_2^-$ - мультипликативные группы алгебр $RP_2^+$, $RP_2^-$, а $R^+ = \{\rho \in R | \rho > 0\}$. Более того, применяя генератор $l_{i,i+1} = diag[1_{(1)},...,l_{(i,i+1)},...,1_{(n)}]_n$ и аналогично построенные генераторы $t_{i,i+1} = diag[t_{(i,i+1)}]_n$, $L_{i,i+1} = diag[L_{(i,i+1)}]_n$, $T_{i,i+1} = diag[T_{(i,i+1)}]_n$, формируем группу четно-составных стрелочных подстановок $n$-ой степени $P_n^\pm = F(\{l_{k,k+1}, \pm t_{s,s+1}\}_I)$ и группу автоморфизмов соответствующей групповой алгебры $Aut\, RP_n^\pm = F(\{L_{k,k+1}, T_{s,s+1}\}_I)$, где $\forall j \in J$ выполняются условия принадлежности $s, k, (k+1) \in I_j$ а $(s+1) \in I_{j+1}$. Таким образом, поскольку конечно-порожденная группа $Aut\, RP_n^\pm$ порядка $n$ образована прямыми суммами $(n-1)$ своих подгрупп 2-го порядка, то, следовательно, если $m=1$, тогда $Aut\, RP_n^\pm = Aut\, RP_n^+ = SO(n)$, если же $m=2$, тогда $Aut\, RP_n^\pm = SO(n_1, n_2)$, что позволяет обобщить понятие специальной ортогональной группы, считая, что $SO(n_1,...,n_m) = Aut\, RP_n^\pm$.

## 2. Многомерные графы и группы топологических автоморфизмов некоторых фактормножеств

Пусть $d$ – произвольное конечное множество одиночных (не расслаиваемых на пары, тройки и т.п.) элементов с отмеченной точкой, тогда $d$ изоморфно $n$-мерному графу, ассоциируемому с простым базисом $R^n$, или подмножеству натуральных чисел с отмеченным нулем, а именно, $d = \{0, I\}$, где $I = \{1,...,n\}$, $n \in N$. Тем самым, имеет место изоморфизм $Aut° d \approx S_n$, где группа автоморфизмов $d$ с закрепленной отмеченной точкой изоморфна симметрической группе, которая, как известно, может быть представлена множеством всех подстановок $n$-ой степени, полученных композицией транспозиций, взятых из набора $\{(i,j)\}_I$. Если же $d$ факторизовано так, что $\forall i, j \in I$ имеем эквивалентность $i \sim j$, тогда, исходя из топологических соображений, 2-цикл $(i,j)$ уже не может быть использован в качестве элементарной подстановки, потому что он предполагает разрыв 2-мерного элемента графа при переходе от графа №1 к №2. В то же время 3-цикл $(i,j,k)$, который является генератором четных подстановок, вполне может быть использован в качестве элементарной подстановки, поскольку его реализация не предполагает разрыва 3-мерного элемента графа №3 при любом циклическом переходе. В свою очередь, если $d$ факторизовано так, что $\forall i, j \in I$ имеем эквивалентность $i \sim j \sim 0$, тогда топологического запрета на 2-цикл $(i,j)$ не последует, т.к. петли графа №4 свободны, т.е. не связаны между собой. Итак, установлено, что если $(d/\sim) \approx \{0,1\}$, то $Aut(d/\sim) \approx S_n^+$, если $(d/\sim) \approx 0$, то $Aut(d/\sim) \approx S_n$, где $Aut(d/\sim)$ - группа автоморфизмов $d$ с закрепленной отмеченной точкой, соответствующих топологии фактормножества $(d/\sim)$, а $S_n^+$ - знакопеременная группа, т.е. множество всех четных подстановок $n$-ой степени, каждая из которых получена композицией 3-циклов, взятых из набора $\{(i,j,k)\}_I$. Важно также отметить, что 2-цикл $(i,j): i \leftrightarrow j$ и 3-цикл $(i,j,k):$ $(i \rightarrow j, j \rightarrow k, k \rightarrow i)$ представляют собой принципиально различные типы отображений (а именно, параллельное и последовательное), которые отмечены обратимой стрелкой и последовательностью однонаправленных стрелок соответственно, в связи с чем, при формировании топологических автоморфизмов фактормножеств действует принцип запрета обратимой стрелки, т.е. принято, что если $i \sim j$ и 0 не эквивалентен им, то элементарная перестановка $i \leftrightarrow j$ не допустима, но возможен переход $i \rightarrow j$ или $j \rightarrow i$ в составе элементарного 3-цикла.

Пусть *b* – произвольное конечное множество двойных (расслаиваемых на пары) элементов вместе с одиночной отмеченной точкой, тогда *b* изоморфно *n*-мерному графу, ассоциируемому с расширенным (простым и обратным к нему) базисом $R^n$, или симметрическому подмножеству целых чисел, т.е. $b \approx \{0, \bar{I}\}$, где $\bar{I} = \{\pm 1,..,\pm n\}$. Тем самым, имеем изоморфизм $Aut^\circ b \approx P_n$, где $P_n$ – группа стрелочных подстановок *n*-ой степени, полученных композицией стрелочных транспозиций, взятых из набора $\{(\pm i), (\pm j, \pm k)\}_I$, где $(\pm i): (\pm i \to \mp i): (+i \leftrightarrow -i)$ и $(\pm j, \pm k): (\pm j \leftrightarrow \pm k)$. Если же *b* факторизовано так, что $\forall i \in I$ имеем эквивалентность $(-i) \sim (+i)$, тогда инверсия не может быть использована в качестве стрелочной транспозиции (действие принципа запрета обратимой стрелки), и поэтому в качестве элементарной стрелочной подстановки служат генераторы четно-инверсных подстановок $(\pm j, \pm k)$, $(\pm j, \mp k)$, что соответствует топологии 3-мерного графа №5. В свою очередь, если *b* факторизовано так, что $\forall j,k \in I$ имеем эквивалентность $(-j) \sim (+j) \sim (-k) \sim (+k)$, тогда выполняется полный запрет обратимой стрелки, и поэтому элементарной подстановкой служит генератор четных стрелочных подстановок 4-цикл $(+j, +k, -j, -k)$, что соответствует топологии 3-мерного графа №6. Если же *b* факторизовано так, что $\forall j,k \in I$ имеем эквивалентность $(-j) \sim (+j) \sim (-k) \sim (+k) \sim 0$, тогда запрета обратимой стрелки нет, и поэтому стрелочная транспозиция не трансформируется. Итак, установлено, что если $(b/\sim) \approx \{0, I\}$, то $Aut(b/\sim) \approx P_n^-$, если $(b/\sim) \approx \{0,1\}$, то $Aut(b/\sim) \approx P_n^+$, если $(b/\sim) \approx \{0, J\}$, то $Aut(b/\sim) \approx P_n^\pm$, если $(b/\sim) \approx 0$, то $Aut(b/\sim) \approx P_n$.

## 3. Поле дискретных чисел и группы автоморфизмов некоторых факторпространств

Пусть $|z|\mod 2: Z \to Z/\sim (z_1 \equiv z_2 (\mod 2) \Leftrightarrow z_1 \sim z_2): |z|\mod 2 = r$, где *r* – остаток от деления *z:2*, взятый без знака *z*, тогда имеем поле $D = \langle \{0,1\}; \oplus, \cdot, 0, 1 \rangle$, где $x \oplus y = |x+y| \mod 2$. Если же *r* – остаток от деления *z:2*, взятый со знаком *z*, тогда имеем не ассоциативную по сложению алгебраическую структуру $B = \langle \{-1,0,1\}; \oplus, \cdot, 0, 1 \rangle$. Пусть $D^n = \prod_n D$, тогда $D^n$ является абелевой группой с покомпонентным сложением. Более того, $D^n$ - линейное пространство над *D* а $e_I = (e_i)_I$, где $e_i = (0_{(1)},...,1_{(i)},...,0_{(n)})$, - его простой базис. Пусть также $H_n^+$ - подгруппа группы $D^n$, состоящая из всевозможных четных сумм элементов простого базиса, а именно, $H_n^+ = F(\{e_j \oplus e_k\}_I)$, где *F(S)* – свободная группа, образованная множеством *S*. Тогда имеем подгруппу $H_n^+ = \{c_n^0,...,c_n^{2m},...,c_n^{2[n/2]}\}$, где $c_n^0 = (0_{(1)},...,0_{(n)})$ а $c_n^{2m} = \{c_I^{2m}(1)\}_{C_n^{2m}}$, где $c_I^{2m}(1)$ - произвольное размещение *2m* единиц и *(n-2m)* нулей, которая разбивает $D^n$ на два класса смежности по этой подгруппе, так что $D^n/H_n^+ \approx D$. Если же принять, что $H_n^\pm = \coprod_J F(\{e_i \oplus e_k\}_{I_j})$, тогда $D^n / H_n^{\pm} \approx D^m$. Одновременно, пусть $B^n = \prod_n B$, тогда $B^n$ - не ассоциативное линейное пространство над $B$ с тем же простым базисом что и $D^n$, но ни $B^n$ ни $B$ не являются группой, и поэтому факторизуются только как множества, т.е. через отношение эквивалентности. Вместе с тем, всякая факторгруппа $D^n / H_n$ изоморфна факторгруппе $Z^n / ZH_n$, где $ZH_n : H_n \to Z^n$ $(0 \to 2Z, 1 \to 2Z+1)$, и тем самым индуцирует соответствующие фактормножества $(Z^n / \sim)$ и $(B^n / \sim)$, причем $(Z^n / \sim) \approx (B^n / \sim)$. Например, если $H_n = H_n^- = 0$, то $ZH_n^- = (2Z)^n$, так что $D^n / H_n^- \approx Z^n / ZH_n^- \approx D^n$, и поэтому $(Z^n / \sim) \approx (B^n / \sim) \approx \{D^n\}$. Таким образом, принимая обозначения, согласованные с принятыми ранее для множеств, изоморфных многомерным графам, а именно, $d = \{0, e_I\}$, где $d \subseteq D^n$, и $b = \{0, \overline{e}_I\}$, $\overline{e}_I = (\pm e_i)_I$, где $b \subseteq B^n$, а также, принимая во внимание, что всякая факторгруппа $D^n / H_n$ индуцирует разбиение $d$ и $b$ на классы эквивалентности, получаем следующие соответствия: если $H_n = H_n^{\pm}$, то $(d / \sim) \approx (b / \sim) \approx \{0, J\}$, а если $H_n = D^n$, тогда $(d / \sim) \approx (b / \sim) \approx 0$. В свою очередь, поскольку автоморфизмы $D^n$ и $B^n$ определяются допустимыми заменами их базисов, то $Aut^\circ D^n = Aut^\circ d = S_n$ и $Aut^\circ B^n = Aut^\circ b = P_n$, а следовательно, определены автоморфизмы соответствующих дискретных факторпространств, а именно: если $H_n = H_n^{\pm}$, то $Aut(B^n / \sim) = P_n^{\pm}$, если же $H_n = D^n$, то $Aut(B^n / \sim) = P_n$. Далее, имея ввиду, что всякому непрерывному факторпространству $R^n / ZH_n$ соответствует некоторое компактное топологическое пространство, симметрии которого определяются факторгруппой $Z^n / ZH_n$, сформулируем утверждение, вытекающее из сравнения классов симметричных точек компактных пространств с фактормножеством $\{Z^n / ZH_n\}$.

**Lemma:** $R^n / ZH_n^+ \approx S^n$, $R^n / ZH_n^{\pm} \approx S^{n_1} \times .. \times S^{n_m}$, $R^n / Z^n \approx RP^n$.

Кроме того, принимая во внимание, что $RP_n \approx End\, R^n$, где $RP_n$ - групповая алгебра над $R$, а $End\, R^n$ - алгебра эндоморфизмов $R^n$, а также то, что расширенный базис $R^n$ изоморфен *n*-мерному графу *b*, сформулируем утверждение об эндоморфизмах факторпространств.

**Lemma:** $RP_n^+ \approx End\, S^n$, $RP_n^{\pm} \approx End\, S^{n_1} \times .. \times S^{n_m}$, $RP_n \approx End\, RP^n$.

Наконец, если учесть, что $Aut\, RP_n^{\pm} = SO(n_1, ..., n_m)$, то тут же получим утверждение об изоморфизме групп.

**Lemma:** $Aut\, S^n \approx SO(n)$, $Aut\, S^l \times S^m \approx SO(l, m)$.

В заключение заметим, что поскольку $P_n^{\pm}$ не исчерпывает всех подгрупп $P_n$, то представленное здесь исследование не полно. Однако, ввиду особого места унитарных пространств в современной физике, дополним его нижеследующим построением. Пусть $C$ – поле комплексных чисел, тогда имеет место овеществление алгебры $C$ в $RP_2^+$, а именно, $\pi: C \to RP_2^+$:

$(x \cdot 1 + y \cdot i) \to \begin{bmatrix} x & y \\ -y & x \end{bmatrix}$, причем $Aut\, C = U(1) \approx Aut\, RP_2^+ = SO(2)$. Кроме того, если поэлементное овеществление матриц обозначить тем же символом, тогда $\pi(u^* u) = \pi(u)^T \pi(u) = diag[1]_{2n}$, где $u \in SU(n)$, а следовательно $\pi(SU(n)) < SO(2n)$. С другой стороны, если принять, что $P_2^\oplus = F(u_1, u_2, u_3, u_4)$, где $u_1 = \begin{bmatrix} 1 & 0 \\ 0 & 1 \end{bmatrix}$, $u_2 = \begin{bmatrix} 0 & i \\ i & 0 \end{bmatrix}$, $u_3 = \begin{bmatrix} 0 & 1 \\ -1 & 0 \end{bmatrix}$, $u_4 = \begin{bmatrix} i & 0 \\ 0 & -i \end{bmatrix}$, а

$RP_2^\oplus = \begin{bmatrix} x+iy & r+is \\ -r+is & x-iy \end{bmatrix}$, $x, y, r, s \in R$, и $U \equiv Aut\, RP_2^\oplus = GP_2^\oplus / R^+ = SU(2)$, где $GP_2^\oplus = U(2)$, тогда, если $U_j = diag[U_{(j,j+1)}]_n$, то $Aut\, RP_n^\oplus = F(\{U_j\}_{I \setminus n}) = SU(n)$,

причем, раскрывая символы $1 = \begin{bmatrix} 1 & 0 \\ 0 & 1 \end{bmatrix}$ и $i = \begin{bmatrix} 0 & 1 \\ -1 & 0 \end{bmatrix}$, получим, что $P_n^\oplus < P_{2n}^+$, где $P_n^\oplus = F(\{(P_2^\oplus)_j\}_{I \setminus n})$ и, как прежде, $(P_2^\oplus)_j = diag[(P_2^\oplus)_{(j,j+1)}]_n$. Таким образом, специальная унитарная и специальная ортогональная группа имеют единую алгебраическую структуру, связанную с автоморфизмами факторпространства $R^{2n} / ZH_{2n}^+$, так что унитарное пространство вкладываются в истинно-евклидово пространство двойной размерности.

***Вывод:*** Жесткость евклидовых пространств, обусловленная определенным типом функционала скалярного произведения, есть результат ограничения, накладываемого на автоморфизмы арифметического пространства при его факторизации. В частности, геометрия истинно-евклидова пространства ассоциируется с факторизацией *n*-мерного арифметического пространства в *n*-мерную сферу.